\documentclass{article}
\usepackage[utf8]{inputenc}
\usepackage[english]{babel}
\usepackage{csquotes}
\usepackage{biblatex}
\usepackage{hyperref}

\usepackage{ dsfont }
\usepackage{amssymb}
\usepackage{amsmath}        
\usepackage{amsfonts}       
\usepackage{amsthm}         
\usepackage{bbding}         
\usepackage{bm}             
\usepackage{graphicx}       
\usepackage{fancyvrb}       
\usepackage{indentfirst}    
\usepackage{icomma}         
\usepackage{dcolumn}        
\usepackage{booktabs}       
\usepackage{paralist}       
\usepackage{xcolor}         
\usepackage{ dsfont }

\input xy                         
\xyoption{all}
\usepackage{amsmath,amscd}
\usepackage{textcomp}
\theoremstyle{plain}
\newtheorem{veta}{Věta}
\newtheorem{Thm}[veta]{Theorem}
\newtheorem{Prop}[veta]{Proposition}
\newtheorem{Ex}[veta]{Example}

\newtheorem{Lemma}[veta]{Lemma}

\theoremstyle{plain}
\newtheorem{Def}[veta]{Definition}
\newtheorem{Remark}{Remark}

\theoremstyle{remark}

\theoremstyle{plain}

\newcommand{\Q}{\mathbb{Q}}
\newcommand{\Z}{\mathbb{Z}}

\newenvironment{dukaz}{
  \par\smallskip\noindent
  \textit{Proof}.
}{

\rightline{$\qedsymbol$}
}

\usepackage{blindtext}
\title{Endomorphism rings of simple modules and block decomposition}
\author{Dominik Krasula
}

\begin{document}
 \maketitle 
\begin{abstract} A left and right noetherian semiperfect ring R is known to be indecomposable if and only if its factor by the second power of Jacobson radical is. This characterisation is used to study simple R-modules in terms of their Ext groups. It is shown that if R is indecomposable, all its simple modules are either finite or have the same infinite cardinality and their endomorphism rings have the same characteristics. The results are further strengthened in the case when R is quasi-Frobenius.
\end{abstract}

\smallskip


\smallskip

\noindent \textbf{Author:} RNDr. Dominik Krasula

Charles University, Faculty of Mathematics and Physics, Prague, Czech Republic

krasula@karlin.mff.cuni.cz

ORCID: 0000-0002-1021-7364

\smallskip 

\noindent This work is a part of project SVV-2023-260721.

\noindent This research was supported by the grant GA ČR 23-05148S from the Czech
Science Foundation.


\section{Introduction}

Let $R$ be an associative ring with 1. By a classical result in ring theory, $R$ can be expressed as a finite direct product of indecomposable rings, called \textit{blocks},  iff 1 can be written as a sum of orthogonal centrally primitive idempotents. This decomposition is then unique up to isomorphism and permutation. 

Such a decomposition always exists if $R$ is a semiperfect ring. Thus, many questions about semiperfect rings can be passed to an indecomposable case. If $R$ is further left and right noetherian, many questions can be passed to a ring $R/J(R)^2$, where $J(R)$ is its Jacobson radical; see [5]. In particular, Ext groups of simple modules are preserved.

The above observations and some classical results in the block theory of artinian rings allow us to study simple modules in terms of \textit{Ext-paths} between them; see Definition \ref{Ext-path}. This viewpoint is closely related to a \textit{quiver} (or scheme) of a ring studied by V. V. Kirichenko, generalising Gabriel quivers of finite-dimensional algebras; see [6]. The formalism of quivers is omitted in this article, but the main results use Kirichenko's work.

\medskip

The considerations in this text are motivated by recent developments in the study of quasi-Frobenius (QF) rings and their relation to MacWilliams rings.  In the late '90s,  J. A. Wood proved in [12] that a finite Frobenius ring is MacWillaims. Due to this and subsequent results, QF rings rose in prominence in algebraic coding theory.  In the last decade, infinite MacWilliams rings and their generalisation became objects of study by module theorists in their own right; see, for example, [1].

M. C. Iovanov solved the problem of characterising left artinian MacWillimas rings in [4]. The key observation was to decompose the ring as a product of a finite ring and a ring with no finite simple modules. This idea goes back to articles [3] and [11].

The particular case for artin algebras was proved in [3]. One of the steps in the proof was showing that over indecomposable left and right artinian rings, either all simple modules are finite or have the same infinite cardinality. In Section \ref{SecMain}, this result is generalised to left and right noetherian semiperfect rings.  We further show that all their endomorphism rings have the same characteristic; see Theorem \ref{Theorem}.

For a semilocal ring, the cardinality of a simple module is strongly linked to the cardinality of its endomorphism ring, as discussed in  \ref{semilocal}, which is a crucial parameter determining the structure of the lattice of submodules of direct powers of simple rings. 

 It was observed in [11, Prop. 1.3], that for a QF ring $R$, the submodule of $Soc(R_R)$ consisting of finite submodules embeds in $R_R/J(R)$ iff   the submodule of $Soc(_RR)$ consisting of finite submodules embeds in $_RR/J(R)$. They coined the term \textit{finitary Frobenius} for such rings. A QF ring is left MacWilliams iff it is a finitary Frobenius ring iff it is a right MacWilliams ring; [11, Cor. 4.8].  The proof of [11, Prop. 1.3] can be easily modified to prove analogous results for the submodule of $Soc(R)$ consisting of simple modules of fixed infinite cardinality. Using the methods of this proof, we show that simple modules linked by a Nakayama permutation have isomorphic rings of endomorphisms; Proposition \ref{PropQF}. 

\medskip

The structure of the article is as follows. Section \ref{SecPrel} gathers known theory on semiperfect and QF rings needed in the sequel. Subsection \ref{SecDecomp} formalises the concept of Ext-path and states and proves Lemma \ref{Lemma} on which most of the main results depend.

Section \ref{SecAbel}  studies Ext groups of modules in terms of characteristics of their endomorphism rings. This characteristic is strongly related to an abelian group structure of a module, which allows one to decide (non)-triviality of the Ext group by studying their extensions in the category of abelian groups. The results of this section are formulated more generally than needed for the paper's main results. The presented method works in a general setting without needing more complicated proofs.

The aforementioned main results of the paper are then discussed in Section \ref{SecMain}.

\section{Preliminaries}\label{SecPrel}

This section gathers the properties and terminology of semiperfect rings needed in Section \ref{SecMain}. We refer the reader to [8, Chapters 7 and 8] for a more detailed exposition and any undefined terminology. We use a different notation for simple modules than in the cited literature. 

The properties following from the semilocality of semiperfect rings are discussed in Subsection \ref{semilocal}. In  Lemma \ref{Lemma}, a classical criterion for indecomposability of radical-square zero artinian rings is reformulated in terms of Ext groups of simple modules. Using the results of [6], Lemma \ref{Lemma} is then generalised to left and right noetherian semiperfect rings; see Theorem \ref{main}. Subsection \ref{SecQF} then gathers the termilogy of Quasi-Frobenius rings.

\medskip 

Throughout this section, $R$ denotes a semiperfect ring, and $J(R)$ is its Jacobson radical. The radical of a finitely generated $R$-module $M$ is denoted by $J(M)$. The factor $M/J(M)$ is called the \textit{top} of $M$. A ring is semiperfect iff any finitely generated module $M$ has a projective cover, denoted by $P(M)$. 

 The operations \textit{top} and \textit{projective cover} induce mutually inverse bijections between isomorphism types of simple modules and isomorphism types of indecomposable f. g. projective modules.  As a module over itself, $R$ decomposes as follows
 \begin{gather}
     R_R\cong P(S_1)^{\mu_1}\oplus \dots \oplus  P(S_n)^{\mu_n} \label{ModuleDecomp} \\ 
          _R R\cong P(_1S)^{\mu_1}\oplus \dots \oplus  P(_nS)^{\mu_n}
 \end{gather}
 where $\mu_i$ are positive integers, $S_1,\dots, S_n$ ($_1S, \dots,~_nS)$ is a list of representatives of isomorphism types of simple right (left) $R$-modules. The decomposition is unique up to isomorphism and permutation. 
 
 Each f.g. indecomposable projective module is a cyclic module generated by a local idempotent. Therefore, decompositions (1) and (2) correspond to a decomposition of 1 as a sum of local pairwise orthogonal idempotents. The formalism of idempotents is omitted in this article, but some results cited from [8] are formulated in this setting. 

\subsection{Semilocal rings}\label{semilocal}

A ring is \textit{semilocal} if its factor by Jacobson radical is semisimple (also called totally decomposable in literature). Over a semisimple ring, simple modules have no nontrivial extensions. For any finitely generated $R$-module $M$, the equality  $J(R)M=J(M)$ holds; [8, Prop. (24.4)].

The Wedderburn-Artin theorem describes the block decomposition of a semisimple ring as a product of finitely many matrix rings. In particular,
\begin{gather}
    R/J(R)\cong M_{\mu_1}(D_1)\times \dots M_{\mu_n}(D_n), \label{TopDecomp}
\end{gather}
where $D_i$ is a division ring isomorphic to $End_R(S_i)\cong End_R(_i S)$.  

 Thus simple right (left) modules are isomorphic to \textit{row} (\textit{column}) modules. That is, $S_i$ ($_iS)$ can be represented as a set of column (row) vectors with $\mu_i$ coordinates with values in $D_i$.  In particular, over a simple semisimple ring, f.g. modules are either all finite if $D$ is, or they have the same infinite cardinality as $D$. This is not true for a general ring, as shown by the following example from [7, Example 8]
\begin{Ex}
     Let $F$ be a finite field and $\kappa$ an infinite cardinal. Then $S:=F^{(\kappa)}$ is an infinite simple  $End_F(S)$-module with a finite ring of endomorphisms.
\end{Ex}

Semiperfect rings are semilocal, so we can view their simple modules as simple modules over the semisimple ring $R/J(R)$. This gives a natural correspondence between left and right simple $R$-modules, as modules $S_i$ and $_iS$ are modules over the same matrix ring in the decomposition of $R/J(R)$. 

As follows from (3), a finite simple module $S_i$ has the same cardinality as its endomorphism ring iff $\mu_i=1$. A semiperfect ring is called \textit{basic} if $\mu_i=1$ for any $i\leq n$.

We end this section with the following lemma used later in the proof of Theorem \ref{Theorem}. It is a generalisation of observation about left and right artinian rings by Iovanov [3, Prop. 2.2]. We present a complete proof for the reader's convenience. 
\begin{Lemma}\label{annihilator} \label{CorCard} Let $R$ be a semilocal ring and  $S$ and $T$ be simple non-isomorphic $R$-modules such that $Ext_R^1(S, T)\neq 0$, and let $I$ and $J$ be their respective annihilator ideals.

Then $IJ$ is strictly contained in $I\cap J$.
\end{Lemma}
\begin{dukaz}
   Because $I$ and $J$ are distinct maximal ideals in $R$, factors $I/IJ$ and $J/IJ$ are distinct maximal ideals in $\bar{R}:=R/IJ$.  Hence $\bar{R}/I\times \bar{R}/J$ is a semisimple ring.    Because $IJ\subseteq I\cap J$, the map  $[a]\mapsto ([a],[a])$ is a monomorphism $\bar{R}/I\cap J\rightarrow \bar{R}/I\times \bar{R}/J\rightarrow 0$.  
   
   If $I\cap J=IJ$, then $\bar{R}$ injects in $\bar{R}/I\times \bar{R}/J$, and hence it is also a semisimple ring.  But by our assumption, $S$ and $T$ have nontrivial extensions as $R$-modules.  Such extensions are annihilated by $IJ$, so we can view them as nontrivial extensions of $\bar{R}$-modules, contradicting that $R/IJ$ is a semisimple ring. 
\end{dukaz}



\subsection{Decomposition of semiperfect noetherian rings} \label{SecDecomp}

This section gathers the basics of the block theory of left and right noetherian semiperfect rings needed to prove the main results. The theory is well known for artin algebras; see [2, Section II.5]. The commutative case is trivial since any commutative semiperfect ring is a finite product of local rings [8, Thm. (23.11)].

If $M$ is an indecomposable $R$-module, there is only one block $B$ whose action on $M$ is nontrivial. In such a case, we say that $M$ \textit{belongs to the block B}. A simple module $S$ belongs to the block $B$ iff it is an epimorphic image of $B$ in $mod\text- R$. Modules that belong to different blocks have only trivial extensions.  

 We say that a simple module $S$ is a \textit{common composition factor} of two modules if, for both modules, there exists a quotient of two consecutive members of a composition series isomorphic to $S$.

\begin{Def} \label{Ext-path}
Let $S$ and $T$ be two simple $R$-modules.

A sequence of simple modules $S=S'_1, S'_2,\dots, S'_k=T$ is called an \emph{Ext-path from $S$ to $T$} if $Ext_R^1(S'_i,S'_{i+1})\neq 0$ or $Ext_R^1(S'_{i+1},S'_i)\neq 0$ for any $i< k$.
\end{Def}
For a basic finite-dimensional algebra $A$ over an algebraically closed field, there is an Ext-path from a simple module $S$ to $T$ iff there is an unoriented path between the corresponding vertices of the Gabriel quiver of $A$ iff they belong to the same block; [2, Prop. III.1.14].

\begin{Lemma}\label{Lemma}
    Let $R$ be a right artinian ring such that $J(R)^2=0$ and $S$, $T$ are two simple right $R$-modules. TFAE

    (1) Modules $S$ and $T$ belong to the same block.

    (2) There is a sequence of indecomposable f.g. projective modules $\mathcal{P}\colon P(S)=P_1, \dots, P_k=P(T)$ such that any two consecutive members have a common composition factor. 

    (3) There is a sequence of indecomposable f.g. projective modules $\mathcal{P}'\colon P(S)=P_1, \dots, P_l=P(T)$ such that their tops induce an Ext-path from $S$ to $T$.

     (4) There is an Ext-path from $S$ to $T$.
\end{Lemma}

\begin{dukaz}
 To show that (3) implies (4), consider the tops of the modules in $\mathcal{P}'$.  Implication (4)$\longrightarrow$(1) can be proved by contraposition since there are no nontrivial extensions of simple modules belonging to different blocks. Implication (1)$\longrightarrow$(2) is a classical result about right artinian rings. See, for example, [8, Thm (22.6)].

Implication (2)$\longrightarrow$ (3): We modify the sequence $\mathcal{P}$. If two consecutive modules are isomorphic, we erase one of them. No change is needed if there is a common composition factor of two consecutive modules isomorphic to the top of one of them. Otherwise, we add a projective cover of a common composition factor between them. We denote the new sequence by $\mathcal{P}'$. 

 Let $P(V)$ and $P$ be two consecutive members of $\mathcal{P}'$. By constructing $\mathcal{P}'$, the simple module $V$ is a composition factor of $P$. Because $P$ and $P(V)$ are non-isomorphic, $V$ is not the simple top of $P$. 

 Module $P$ is finitely generated, so $0=J(R)^2P=J^2(P)$, and hence,  $J(P)$ is a semisimple module. In particular, there is a projection $p\colon J(P)\twoheadrightarrow V$. The canonical projection $P(V)\twoheadrightarrow V$ lifts to an epimorphism $P(V)\twoheadrightarrow J(P)$ that commutes with $p$. As a result, $P/Ker(p)$ is a nontrivial extension of $Top(P)$ by $V$.
\end{dukaz}

\begin{Thm} \label{main}
    Let $R$ be a semiperfect ring and further assume that one of the following conditions holds

    (1) $R$ is left and right noetherian.

    (2) $R$ is left and right perfect and $R/J^2(R)$ is left and right noetherian.
    
    Then $R$ is indecomposable iff for any two simple right $R$-modules $S$ and $T$ there is an Ext-path from $S$ to $T$. 
\end{Thm}
\begin{dukaz}
By  [6, Thm. 1.2], left and right Noetherian ring $R$ is indecomposable iff $R/J(R)^2$ is. The same holds for a left and right perfect ring [6, Thm. 3.4].

 By the correspondence theorem for rings, $J(R)/J(R)^2$ is the Jacobson radical of $R/J(R)^2$. So $R/J(R)^2$ is a semilocal ring, and its radical is nilpotent. If $R$ is right and left noetherian, $R/J(R)^2$ is right and left noetherian, hence right and left artinian by the Hopkins-Levitzki theorem. 

    Recall that $J(R)$ annihilates all simple modules, and so $J(R)^2$ annihilates any extension of a simple module by a simple module. Therefore, we can identify extensions of simple $R$-modules with extensions of simple $R/J(R)^2$-modules.  Lemma \ref{Lemma} then completes the proof. 
\end{dukaz}
We end with a couple of examples showing the limits of the above theorem and, hence the limits of its corollaries from Section \ref{SecMain}
\begin{Ex}
    Consider a commutative noetherian ring  $\Z$. It is not semilocal and hence not semiperfect. 

   The ring of integers is indecomposable, i.e., all its simple modules are in the same block. But the extension group of any two simple non-isomorphic modules, $\Z_p$ and $\Z_q$, is trivial by a standard calculation. 
\end{Ex}
\begin{Ex}[Kirichenko]\label{KirichenkoEx}
   For a fixed prime $p$, consider a ring $A=\begin{pmatrix}
        \Z_{(p)} & \Q\\
        0 & \Q 
    \end{pmatrix}$, with radical $J(A)=\begin{pmatrix}
        p\Z_{(p)} & \Q\\
        0 & 0
    \end{pmatrix}$, where $\Z_{(p)}$ is the localisation of $\Z$ at $p$.

    Ring $A$ is semiperfect right but not left noetherian. It is indecomposable as a ring, but $A/J(A)^2\cong \Z_p^2\times \Q$ is decomposable.     Two simple modules, $E_{1,1}A\cong \Z_p $ and $E_{2,2}A\cong\Q$, lie in the same block, but their extensions split.     
\end{Ex}
\subsection{QF rings}\label{SecQF}

A ring $R$ is quasi-Frobenius, or QF for short,  if  $D:=Hom_R(-, R)$, given by the $R\text- R$ bimodule $R$, is a duality between the left and right categories of finitely generated $R$-modules. Such a ring is necessarily a two-sided artinian. All projective modules are injective, and vice versa. In particular, indecomposable f.g. projectives have simple socles. See [9, Chapter 6] for details. 

A left and right artinian ring is QF iff there exists a \textit{Nakayama permutation}, i.e., a permutation of $\{1,\dots, n\}$ such that $Soc(P(S_i))\cong S_{\pi(i)}$ and $Soc(P(_{\pi(i)}S))\cong {_iS}$.  We then say that simple modules $S_i$ and $Soc(P(S_i))$ (and their left counterparts) are \textit{linked} by the Nakayama permutation.

\section{Group structure of a module and Ext groups} \label{SecAbel}

This section shows that the triviality of the group of extensions of two modules can often be deduced by viewing the module as an abelian group, i.e., \textit{forgetting} the structure given by an action of a ring and viewing them in $Mod\text- \Z$.  Some group properties of a module, such as divisibility and torsion can be deduced from the characteristics of its ring of $R$-endomorphisms. 

Let $R$ be an arbitrary ring and $M$ a left $R$-module. By  \textit{the characteristic} of $M$, denoted by $char(M)$, we mean the characteristic of the ring $End_R(M)$.  If  $char(M)=a>0$, then $aM=0$. 

When we want to stress that we view an $R$-module as an abelian group, we use the subscript $_\Z$. We refer the reader to [10] as a reference for group properties used in the proofs.

\begin{Prop} \label{Characteristic}
Let $R$ be a ring and $M$, $N$ two left $R$-modules and $a,b$ integers such that $char(M)=a$ and $char(N)=b$.

(1) If  $a$, $b$ are positive and coprime, then $Ext_R^1(N,M)=0$.

(2) If $M$ is simple and $a=0\neq b$ and $M$, then  $Ext_R^1(N,M)=0$.

(3) If  $M$ and $N$ are simple modules and $a\neq 0=b$, then $Ext_R^1(N,M)=0$.
    
\end{Prop}

The method of proof is as follows: we consider an $R$-module extension
\begin{gather}
    \epsilon\colon ~ 0 \rightarrow M \xrightarrow{\nu} K \xrightarrow{\pi} N\rightarrow 0.  \label{epsilon}
\end{gather}
The characteristics of modules $M$ and $N$ give us information about the orders of their elements, which we use to prove that  $\epsilon$ splits in $Mod\text-\Z$, i.e., there is a group decomposition 
\begin{gather}
    K\cong \nu(M) \oplus \pi'(N) \label{split}
\end{gather}
where $\pi'$ is a one-sided inverse of $\pi$ in the category of abelian groups.  We prove that $\pi'(N)$ is a fully characteristic subgroup of $K$, i.e., invariant under any group endomorphism of $K$. In particular, $\pi'(N)$ is invariant under all endomorphisms of the form $r\cdot$ for $r\in R$. Thus $\pi'(N)$ is an $R$-submodule and hence direct summand in $Mod\text-R$, proving that $\epsilon$ splits in $Mod\text- R$.

\begin{dukaz}
    
\textbf{Proof of part (1)} By the Prüfer-Baer theorem, any abelian group of bounded order is a direct sum of cyclic groups; [10, Cor. 10.37]. The orders of all elements of $_\Z N$ are bounded by $b$. Therefore, $Ext^1_\Z(N, M)$ is isomorphic to a direct product of groups of the form $Ext^1_\Z(\Z_p^k, M)$, where $k$ is a positive integer and $p$ is a power of a prime number such that  $p\mid b$.  By our assumptions on $a$ and $b$, the order of any element of the group $_\Z M$ is coprime with $p$. By standard calculation, $Ext^1_\Z(\Z_p^k,M)=0$ and hence any extension of $N$ by $M$ splits in $Mod\text-\Z$.

Take  $\epsilon$  as in (4) and (5).  Because maps $\nu$ and $\pi'$  are injective, they preserve orders of elements. In particular, $K$ is a torsion group, so it decomposes as a direct sum of $p$-primary groups [10, Thm. 10.7]. For any $m\in \nu(M)$ and $n\in \pi'(M)$, the orders of elements $m$ and $n$ are coprime by the assumptions on $char(M)$ and $char(N)$. Hence, each lies in a different $p$-primary component.  We conclude that $\pi'(N)$ is a direct sum of some $p$-primary components of $K$, and thus it is a fully characteristic subgroup.

\textbf{Proof of part (2)} Because $char(M)=0$, the map $nId_M$ is nonzero for any positive $n$, and hence an automorphism by the simplicity of $M$. The surjectivity of maps $nId_M$ implies that $_\Z M$ is divisible; their injectivity implies that $_\Z M$ is torsion-free. 

Consider  $\epsilon$  as in (4). Because $_\Z M$ is divisible, the extension $\epsilon$ splits in $Mod\text-\Z$, as in (\ref{split}). The subgroup $\pi'(N)$ is the torsion subgroup of $K$ and hence a fully characteristic subgroup.

\textbf{Proof of part (3)} Consider  $\epsilon$  as in (4). Because $N$ is a simple module and $char(N)=0$, the group $K/\nu(M)\cong N$ is divisible. That is, for any positive $q$, we have $q(K/\nu(M))=K/\nu(M)$, giving an $R$-module decomposition $qK+\nu(M)=K$. We show that $qK\cap \nu(M)=0$ for some $q$ and hence $\epsilon$ splits.

If $Kq=K$ for all positive $q$, then $_\Z K$ is divisible, and its torsion subgroup $_\Z \nu(M)$  is divisible too. But $aM=0$, which is a contradiction. Thus, we can assume that there is $q>0$ such that $qK$ is a proper simple submodule of $K$. Group $_\Z K$ is not a torsion group, so $aqK\neq 0$ and $qK\neq \nu(M)$. Because both $\nu(M)$ and $qK$ are simple modules, $qK\cap \nu(M)=0$, and this finishes the proof.  

\end{dukaz}

\begin{Remark}
The assumption of simplicity in parts (2) and (3) is not superfluous. By a standard calculation, $Ext_\Z^1(\Z_m, \Z)\cong\Z_m\neq 0$ for any positive integer $m>1$.  

If the endomorphism rings of modules $N$ and $M$ induce a bimodule structure on their Ext groups, for example, if $R$ is a finite-dimensional algebra over a field, one can omit the assumption that $M$ is simple in part (3).    
\end{Remark}

\section{Main results}\label{SecMain}

This section gathers and proves the main results of this article. 

\begin{Thm}\label{Theorem} Let $R$ be a semiperfect left and right noetherian ring. 

Then there exists a decomposition of $R$ as a finite direct product of blocks such that for each block $B$, the following holds.

(1) If $B$ has a finite simple module, then all finite-length $B$-modules are finite.  All simple modules have equal cardinality as their endomorphism rings iff $B$ is basic. 

(2) If $B$  has no finite simple modules, then there exists infinite cardinal $\kappa$ such that all finite-length $B$-modules have cardinality $\kappa$. Furthermore, endomorphism rings of simple modules have cardinality $\kappa$.

(3) All endomorphism rings of simple modules have the same characteristic.  
\end{Thm}
\begin{dukaz}
Direct summands of semiperfect and noetherian rings are again semiperfect and noetherian, so each block $B$ satisfies the assumptions of Theorem \ref{main}.  Part (3) then follows  Proposition \ref{Characteristic}. 

    Block $B$ is semilocal, so the relations between the cardinality of simple modules and their endomorphism rings follow from Section \ref{semilocal}.

Let $\kappa$ be the largest cardinality of a simple module over $B$. Then any module of composition length $l$ has cardinality at most $\kappa^l$. Thus, to show (1) and (2),  it is enough to prove that all simple modules over $B$ are either finite or have the same infinite cardinality. 

 Let $S$ and $T$ be simple non-isomorphic $B$-modules such that $Ext_B^1(S, T)\neq 0$, and let $I$ and $J$ be their respective annihilator ideals. By Lemma \ref{annihilator}, $IJ$ is strictly contained in $I\cap J$.  Because $B$ is left and right noetherian, $I\cap J$ is finitely generated as both right and left $B$-module. In particular, $I\cap J/IJ$ is finitely generated as both an $B/I$-module and an $B/J$-module. Both rings $B/I$ and $B/J$ are simple semisimple rings. Therefore, $B/I$, $I\cap J/IJ$, and $B/J$ are either all finite, and so $S$, $T$ are
finite, or have the same (infinite) cardinality, equal to the cardinality of $S$ and $T$.    
\end{dukaz}
\begin{Remark}
    By Theorem \ref{main},  part (3) also holds in the case when blocks are left and right perfect and factors by the second power of their Jacobson radical is left and right noetherian.
\end{Remark}
Under the assumptions of Theorem \ref{Theorem}, there is an easy criterion to decide whether a finite-length module has finitely many submodules. We decompose it as a direct sum of finitely many indecomposable modules. For each indecomposable module $M$, there is only one block $B$ whose action on $M$ is nontrivial. If the indecomposable module is finite, the conclusion is trivial. In the infinite case, the following proposition applies.
\begin{Prop}
    Let $R$ be an indecomposable semiperfect left and right noetherian ring and $M$ be an infinite module of finite length. 

    Then $M$ has finitely many submodules iff all epimorphic images of $M$ have simple socles.
\end{Prop}
\begin{dukaz}
    By Theorem \ref{Theorem}, endomorphism rings of simple $R$-modules are all infinite. The conclusion then follows from [7, Prop. 14]
\end{dukaz}

In the case of QF rings, we can strengthen Theorem \ref{Theorem}.
\begin{Prop} \label{PropQF}
Let $R$ be a QF ring and $\pi$ a corresponding Nakayama permutation. 

Then $End_R(S_i)\cong End_R(S_{\pi(i)})$ and $End_R(_{\pi(i)}S)\cong End_R(_iS)$. 
\end{Prop}
\begin{dukaz}
    We only prove the version for the right modules; the left version is analogous. 

    Consider the right $R$-module $P(S_{\pi^{-1}(i)})$, then duality $D=Hom_R(-,R)$ maps its simple socle $S_i$  to a simple top $_{\pi^{-1}(i)}S$ of the left $R$-module $P(_{\pi^{-1}(i)}S)$. By Schur's lemma, $D(S_i)\cong Hom_R(S_i, R)\cong Hom_R(S_i, S_i^{\mu_i})\cong End_R(S_i)^\mu_i$, where $\mu_i$ is the multiplicity of $P(S_i)$ in the decomposition of $R_R$, as in (1). Thus viewed as $R/J(R)$-module, $D(S_i)$  is the simple column module with entries in $D_i=End_R(S_i)$. 

    Thus, as $R/J(R)$-module, $_{\pi^{-1}(i)}S$ corresponds to a column module $End_R(S_i)^\mu_i$, so, in particular, they have the same ring of endomorphisms. In total, we get 
    \[End_R(S_i)\cong End_R(_{\pi^{-1}(i)}S)\cong  End_R(S_{\pi^{-1}(i)}),\]
    where the last isomorphism follows from Subsection \ref{semilocal}.
    
\end{dukaz}


\medskip

\noindent\textbf{Aknowledgements:}
 The author wishes to express his thanks and appreciation to Jan Žemlička for his guidance and valuable insights during the research. 

\medskip

\noindent \textbf{Bibliography} 

 [1] P. A. G. Asensio, A. K. Srivastava,  MacWilliams extending conditions and quasi-Frobenius rings, \textit{J. of Algebra}  \textbf{605} (2022) 394-402. 
https://doi.org/10.1016/j.jalgebra.2022.05.005

[2]  M. Auslander,  I.  Reiten, S. Smal\o, \textit{Representation Theory of Artin Algebras} (Cambridge Univ. Press, 1995)

[3] M. C. Iovanov,   Frobenius–Artin algebras and infinite linear codes. \textit{J. Pure Appl. Algebra}. \textbf{220(1)}. (2016). 560-576 https://doi.org/10.1016/j.jpaa.2015.05.030

[4]  M.  C. Iovanov,   On infinite MacWilliams rings and minimal injectivity conditions. \textit{Proceedings of the AMS}. \textbf{150(11)}.  (2022). 4575-4586. https://doi.org/10.1090/proc/15929

 [5] V. V. Kirichenko, Generalized uniserial rings,  \textit{Mat. Sb.} \textbf{99(141)} (1976) 559–581.  https://doi.org/10.1070/SM1976v028n04ABEH001666
 
[6] V. V. Kirichenko, Quivers of associative rings, \textit{J. Math. Sci.} \textbf{131(6)} (2005) 6032-6051 https://doi.org/10.1007/s10958-005-0459-6

[7] D. Krasula, Möbius function for modules and thin representations, \textit{ArXiv}, https://doi.org/10.48550/arXiv.2403.05656

[8] T. Y. Lam, A First Course in Noncommutative Rings (Springer-Verlang, 1991)

[9] T. Y. Lam, Lectures on Modules and Rings (Springer-Verlang, 1998)

[10] J. J. Rotman, An Introduction to the Theory of Groups (Springer-Verlang, 1999

[11]   F. M. Schneider, J. Zumbragel,  MacWilliams’ extension theorem for infinite rings, \textit{Proceedings of the AMS} \textbf{147(3)} (2019) 947-961. 
https://www.ams.org/journals/proc/2019-147-03/S0002-9939-2018-14343-9

[12] J. A. Wood,  Duality for modules over finite rings and applications to coding theory, \textit{Am. J. Math.} \textbf{121(3)} (1999) 555–575. http://www.jstor.org/stable/25098937

\end{document}